\documentclass[a4paper]{amsart}

\usepackage{a4}

\newtheorem{theorem}{Theorem}[section]
\newtheorem{proposition}[theorem]{Proposition}

\theoremstyle{remark}
\newtheorem{remark}[theorem]{Remark}

\newcommand{\bee}[1]{\begin{equation}\label{#1}}
\newcommand{\beq}[1]{\begin{eqnarray}\label{#1}}
\newcommand{\ene}{\end{equation}}
\newcommand{\eqe}{\end{eqnarray}}

\begin{document}

\title[Structure theory of PI-rings]{An alternative approach to the structure theory of PI-rings}

\thanks{2010 {\em Math. Subj. Class.} 16R20,  16R50, 16R60.}
\thanks{Supported by the Slovenian Research Agency (program No. P1-0288).
}
\author{Matej Bre\v sar}
\address{Faculty of Mathematics and Physics, University of Ljubljana, and Faculty of Natural Sciences and Mathematics,
 University of Maribor, Slovenia	}
\email{matej.bresar@fmf.uni-lj.si}
\date{}

\begin{abstract}
We expose a rather simple and direct approach to the structure theory of prime PI-rings (``Posner's theorem"), based on fundamental properties of the extended centroid of a prime ring.
\end{abstract}

\maketitle

\maketitle

\section{Introduction}
The theory of rings with polynomial identities originated in Kaplansky's 1948 paper \cite{Kap}, in which he showed that a primitive PI-algebra is finite-dimensional over its center. In 1960 Posner  \cite{Pos} extended  this theorem to the prime ring context; he proved that a prime PI-ring has a two-sided classical ring of quotients which is a finite-dimensional central simple algebra. After the discovery of central polynomials on matrix algebras in the early 1970's, Posner's theorem was further improved by noticing (by different authors, cf. \cite{Row}) that this ring of quotients is actually the algebra of central quotients.

A standard proof of this sharpened version of Posner's theorem, which can be found in several graduate algebra textbooks  (e.g., in \cite{Bah, MR, Row2}), is a beautiful illustration of the power and applicability of the classical structure theory of rings. Its main ingredients are the Jacobson density theorem, the theorem by Nakayama and Azumaya on maximal subfields of division algebras, Amitsur's theorem on the Jacobson radical of the polynomial ring, the nonexistence of nonzero nil ideals in semiprime PI-rings, and the existence of  central polynomials on matrices. The first  two theorems are  needed for the proof of Kaplansky's theorem, which is an  intermediate step in this standard proof of Posner's theorem.

The purpose of this paper is to give a more streamlined proof, which in  each of its steps avoids representing elements in our rings as matrices or linear operators. All aforementioned ingredients are replaced by a single theorem (Theorem \ref{LJerry})
   describing one of the basic properties of the extended centroid of a prime ring. This theorem is essentially due to Martindale \cite{Mart}, and is one of the cornerstones in the theory  of generalized polynomial identities \cite{BMMb} as well as in the theory of functional identities \cite{FIbook}. Our proof is in fact more typical for these two theories than for the PI theory. We have to point out, however, that the idea to use such an approach is not entirely new. Already in \cite{Mart} Martindale noticed that Posner's theorem can be derived from his result on generalized polynomial identities in prime rings (see also \cite{BMMb}). But the proof of the latter is not so easy. Focusing only on polynomial identities, but hiddenly regarding them as generalized polynomial and functional identities, we will be able to get a rather simple and straightforward proof.

Section 2 surveys the prerequisites needed for our proof. In Section 3 we study identities in central simple algebras, and in particular prove Kaplansky's theorem for them. This weak version of Kaplansky's theorem is our   intermediate step which, as we show in Section 4, quickly yields the final result.

\section{Preliminaries}

The purpose of this section is to make this paper accessible to non-specialists.
It is divided into two parts. In the first one we review the  properties of the extended centroid and related notions, and in the second one we give an elementary  introduction  
to polynomial identities. 

\subsection{The extended centroid} Let $R$ be a {\em prime ring}, i.e., a ring in which the product of two nonzero ideals is always nonzero. Then one can construct the {\em  symmetric Martindale ring of quotients} $Q=Q_s(R)$ of $R$, which is, up to isomorphism, characterized by the following four properties: 
\begin{itemize} \item[(a)] $R$ is a subring of $Q$;
 \item[(b)]  for every $q\in Q$ there exists a nonzero ideal $I$ of $R$ such that  $q I\cup I q\subseteq R$;
  \item[(c)] if $I$ is a nonzero ideal  of $R$ and $0\ne  q \in Q$, then $qI\ne 0$ and $Iq\ne 0$;
  \item[(d)] if $I$ is a  nonzero ideal  of $R$,  $f:I\to R$ is a right $R$-module homomorphism, and  $g:I\to R$ is a left $R$-module homomorphism such that $xf(y) = g(x)y$ for all $x,y\in I$, then there exists $q\in Q$ such that $f(y)=qy$ and $g(x)= xq$ for all $x,y\in I$. 
\end{itemize}
 The center $C$ of $Q$ is called the {\em extended centroid} of $R$.   It is a field containing the center $Z$ of $R$. We remark that $Z$ has no zero divisors, and therefore, provided  it is nonzero, one can form its field of fractions. This is  a  subfield of $C$; examples where it is a proper subfield can be easily constructed. 
We may consider $Q$ as an algebra over $C$. A subalgebra of special importance is the so-called {\em central closure} of $R$, which we denote by $R_C$. It consists of  elements of the form $\sum_i \lambda_i r_i$, where $\lambda_i\in C$ and $r_i\in R$. Both $Q$ and $R_C$ are prime rings. The extended centroid of $R_C$ is just $C$. The same is true for every nonzero ideal of $R_C$ (as well as of $R$). If $C\subseteq R_C$, then $C$ is the center of $R_C$. 

The main property of $C$ that we need is given in the following theorem. Its original version was proved by Martindale  in \cite{Mart}. The version that we state is, as one can see from \cite[Theorem A.4]{FIbook}, a special case of
\cite[Theorem A.7]{FIbook}.

\begin{theorem}\label{LJerry}
Let $R$ be a prime ring with extended centroid $C$,  and let $I$ be a nonzero ideal of $R$.
 Assume that $a_i,b_i,c_j,d_j\in Q_s(R)$ satisfy $\sum_{i=1}^n a_ixb_i = \sum_{j=1}^m c_jxd_j $ for all $x\in I$. If $a_1,\ldots, a_n$ are linearly independent over  $C$, then each $b_i$ is a linear combination of $d_1,\ldots, d_m$.
\end{theorem}

Proving Theorem \ref{LJerry}, as well as all other facts about $C$ and $Q$ mentioned above, is neither difficult nor lengthy; it is also entirely self-contained, and can  be easily incorporated into a course on noncommutative rings. See \cite[Chapter 2]{BMMb} for a detailed, and \cite[Appendix A]{FIbook} for a more informal survey on this subject.

If $R$ is a simple ring with $1$, then it follows easily from (a)-(c) that $R=Q$ and hence $C$ is the center of $R$. We may regard every such ring as a central simple algebra (recall that an algebra  over a field is said to be {\em central} if  its center consists of scalar multiples of $1$). Our central simple algebras may be infinite-dimensional.

\subsection{Polynomial identities} Let $C$ be a field, and let $C\langle X_1, X_2,\ldots\rangle$ be the free algebra over $C$ generated by the indeterminates $X_i$, $i=1,2,\ldots$. One can view elements in $C\langle X_1, X_2,\ldots\rangle$ as polynomials in noncommuting indeterminates $X_i$. The {\em degree} of such a  polynomial is defined in a self-explanatory manner.  Let $A$ be an algebra over $C$ and let $f=f(X_1,\ldots,X_n)\in C\langle X_1, X_2,\ldots\rangle$.  We say that $f$ is an {\em identity} of $A$ if $f(a_1,\ldots,a_n)=0$ for all $a_i\in A$. If $f\ne 0$, then $f$ is called a {\em polynomial identity} of $A$. We say that $f$ is a PI-{\em algebra} if there exists a polynomial identity of $A$. 

 An element in $ C\langle X_1, X_2,\ldots\rangle$ of the form
 \begin{equation} \nonumber
 \sum_{\pi \in S_m} \alpha_\pi X_{\pi(1)}\ldots X_{\pi(m)},\quad \alpha_\pi\in C,
 \end{equation}
 where $S_m$ is the symmetric group of degree $m$,
 is called a {\em multilinear polynomial}. Especially important examples are the so-called {\em standard polynomials} in which  $\alpha_\pi$ is defined as the sign of the permutation $\pi$. The standard polynomial of degree $m$ will be denoted by $St_m$.
 The simplest example of a polynomial that is not multilinear is $X_1^2$. However, if  this polynomial is an identity of $A$, then so is the multilinear polynomial $X_1X_2 + X_2X_1= (X_1+X_2)^2-X_1^2 - X_2^2$. Somewhat more tedious, but  based on the same simple idea, is to prove that if $A$ satisfies a polynomial identity of degree $n$, 
 then it also satisfies  a  
 multilinear polynomial identity of degree $\le n$. Accordingly, if we are interested only in the structural properties of a PI-algebra, we may confine ourselves to the study of multilinear polynomials.

A commutative algebra satisfies the polynomial identity $St_2=[X_1,X_2]$. Next, every finite-dimensional algebra is a  PI-algebra. Namely,  if $\dim_C A = m$ then $A$ satisfies $St_{m+1}$. Another important class of PI-algebras is provided by the  Amitsur-Levitzki theorem: If $R$ is any (possibly 
 infinite-dimensional) commutative algebra, then $St_{2n}$ is a polynomial identity of  the matrix algebra $A=M_n(R)$.

Now let $R$ be ``merely" a ring. One can then consider identities of $R$ as elements in $\mathbb Z\langle X_1, X_2,\ldots\rangle$. However, some care is needed in defining when $R$ is a PI-ring. Some trivial polynomials, such as $pX_1$ if $R$ has characteristic $p$, must be excluded. Since we will be interested only in prime rings, we give just the definition adjusted to this context: a {\em prime ring} $R$ is said to be a PI-{\em ring}  if a nonzero polynomial in $C\langle X_1, X_2,\ldots\rangle$, where $C$ is the extended centroid of $R$, is an identity of $R$. 
An illustrative example 
 is $R=M_n(\mathbb Z)$. It satisfies $St_{2n}$, so  is a prime PI-ring. Its extended centroid  is isomorphic to $\mathbb Q$, and its central closure   is isomorphic to $M_n(\mathbb Q)$.

Everything said so far about polynomial identities is the most standard material that can be found in numerous textbooks.

\section{Central simple PI-algebras}

The goal of this section is to prove a proposition that combines two well known results: Kaplansky's theorem on primitive PI-algebras \cite[Theorem 1]{Kap} and Martindale's theorem on prime GPI-rings \cite[Theorem 2]{Mart}. However, we consider only simple algebras in our proposition. The novelty is  a simple proof, adjusted to this special setting.

Let $A$ be an algebra. For  $a,b\in A$ we define $L_a,R_b:A\to A$ by $L_a(x) = ax$, $R_b(x)=xb$. Obviously, $L_aR_b = R_b L_a$. The set $M(A)$ of all operators of the form $\sum_i L_{a_i}R_{b_i}$, $a_i,b_i\in A$, forms a subalgebra of  the algebra  End$_C(A)$ of all linear operators on $A$. We call $M(A)$  the {\em multiplication algebra of $A$}. We remark that Theorem \ref{LJerry} considers two elements in $M(Q_s(R_C))$.

\begin{proposition}\label{P1}
Let $A$ be a  central simple algebra over a field $C$. The following conditions are equivalent:
\begin{itemize}
\item[(i)] $A$ is a PI-algebra;

\item[(ii)] $M(A)$ contains a nonzero finite rank operator;

\item[(iii)] $\dim_C A<  \infty$;

\item[(iv)] $M(A) =  {\rm End}_C(A)$.
\end{itemize}  
\end{proposition}

\begin{proof}(i)$\Longrightarrow$(ii).  Let $f = f(X_1,\ldots,X_n)$ be a multilinear polynomial  
identity of $A$. Pick $1\le i < j\le n$, and write $f=f_i + f_j$ where $f_i$ is the sum of all monomials of $f$ of the form $mX_i m' X_jm''$, and $f_j$ is the sum of all monomials of $f$ of the form $nX_j n' X_in''$
(here, of course, $m,m'$ etc. are monomials in the other variables). Suppose that both $f_i$ and $f_j$ are identities of $A$. Since $f\ne 0$, we have $f_i\ne 0$ or $f_j\ne 0$. Without loss of generality we may assume that $f_i\ne 0$.  Now 
 we may replace the role of $f$ by $f_i$, and hence assume that $X_i$ appears before $X_j$ in all monomials of $f$.   Since $X_1X_2\ldots X_n$ is not an identity of $A$, there exists a pair  $1\le i < j\le n$ such that $f_i$ and $f_j$ are not identities of $A$. We may assume that $i=1$ and $j=2$.

Fix $u_i\in A$ such that $f_1(u_1,\ldots,u_n)\ne 0$.  
The identity $f_1(x,y,u_3,\ldots,u_{n}) = -f_2(x,y,u_3,\ldots,u_{n})$ for all $x,y \in A$ can be rewritten as  a (functional) identity
\begin{equation} \label{eaf}
\sum_{i=1}^r a_i x T_i(y) = \sum_{j=1}^s S_j(y)xd_j\quad\mbox{ for some $T_i,S_j\in M(A)$, $a_i,d_j\in A$.} 
\end{equation}
Since both sides of \eqref{eaf} are nonzero if we take $x=u_1$ and $y=u_2$, some of the $a_i$'s are nonzero. Without loss of generality we may assume that $\{a_1,\ldots,a_t\}$, $t\le r$, is a maximal linearly independent subset of  $\{a_1,\ldots,a_r\}$. Expressing the $a_i$'s with $i >t$ through the $a_i$'s with $i \le t$ we see that  \eqref{eaf} can be rewritten as
\begin{equation} \label{eaf2}
\sum_{i=1}^t a_i x W_i(y) = \sum_{j=1}^s S_j(y)xd_j\quad\mbox{ for some $W_i,S_j\in M(A)$, $a_i,d_j\in A$.} 
\end{equation}
Of course, some of the $W_i$'s are nonzero;  we may assume that $W_1$ is one of them. Now, for any fixed $y\in A$ we infer from \eqref{eaf2} and Theorem \ref{LJerry} that $W_1(y)$ lies in the linear span of $d_1,\ldots,d_s$. Thus, (ii) holds.

(ii)$\Longrightarrow$(iii). Let $W = \sum_{i=1}^n L_{a_i}R_{b_i}$ be a nonzero finite rank operator in $M(A)$.
Picking a maximal linearly independent subset of $\{a_1,\ldots,a_n\}$ and then expressing the other $a_i$'s as linear combinations of elements from this set, we see that there is no loss of generality  in assuming the linear independence of  $\{a_1,\ldots,a_n\}$. We may also assume that $b_1\ne 0$. 
The proof is by induction on $n$.

Let $n=1$. Since  $A$ is simple, there exist $u_j,v_j, w_k,z_k\in A$ such that $\sum_j u_ja_1 v_j =\sum_k w_k b_1 z_k =1$. Consequently, 
$\sum_{j,k} L_{u_j}R_{z_k}WL_{v_j}R_{w_k}$ is the identity operator, and is  of finite rank. But this means that $\dim_C A<  \infty$.

Now let $n > 1$. We will just repeat the appropriate argument from \cite{Mart}. If each $b_i$, $i\ge 2$, is a scalar multiple of $b_1$, then we are back to the $n=1$ case. We may therefore assume that $b_2$ and $b_1$ are linearly independent. By Theorem \ref{LJerry}
there exists $c\in A$ such that $b_1 c b_2 \ne  b_2 c b_1$. Define $W'\in M(A)$ by $W' =  W R_{b_1 c} -R_{cb_1} W$. Obviously, $W'$ has finite rank, and note that 
$
W' = \sum_{i=2}^n L_{a_i} R_{c_i}
$
where $c_i= b_1cb_i - b_icb_1$. Since $a_2,\ldots, a_n$ are linearly independent and $c_2\ne 0$, Theorem \ref{LJerry} shows that $W'\ne 0$.  By induction the proof is complete.

(iii)$\Longrightarrow$(iv). Let $\{a_1,\ldots,a_n\}$ be a basis of $A$. Suppose $\lambda_{ij} \in C$ are such that $\sum_{i,j=1}^n \lambda_{ij} L_{a_i}R_{a_j} =0$. Rewriting this as 
 $\sum_{i=1}^n  L_{a_i} \bigl(\sum_{j=1}^n  \lambda_{ij} R_{a_j}\bigr) =0$ 
 we see by using Theorem \ref{LJerry}  that $\sum_{j=1}^n \lambda_{ij} R_{a_j} =0$, which in turn yields $\lambda_{ij} =0$ for all $i,j$. Therefore $\dim_C M(A) = n^2 = \dim_C {\rm End}_C(A)$, and so  $M(A) =  {\rm End}_C(A)$.
 
 Since (iii)$\Longrightarrow$(i) and (iv)$\Longrightarrow$(ii) are trivial, this completes the proof.
\end{proof}

\begin{remark}\label{R1} 
The first step in the standard proof of Kaplansky's theorem is the reduction to the case where the algebra in question is a division algebra. This can be done quite easily by applying the Jacobson density theorem. Proposition \ref{P1} of course covers division algebras, so we now have  a new proof of Kaplansky's  theorem that does not use the  theorem on maximal subfields of division algebras.
\end{remark}

\begin{remark}\label{R2} 
The proof of (i)$\Longrightarrow$(ii) is also applicable to generalized polynomial identities. Explaining this in detail would make this paper, which is intended for a wider audience, too technical.  Therefore we will just give a few comments that should be sufficient for specialists. Let  $f$ be a multilinear generalized polynomial identity of degree $n\ge 2$ of a prime ring $R$ (in the sense of \cite{BMMb}).  Fixing all variables except two, we arrive at \eqref{eaf}. The only problem is to show that we can choose the two non-fixed variables  in such a way that both sides of  \eqref{eaf} are not identically zero. We argue as in the first paragraph of the (i)$\Longrightarrow$(ii) proof, and in that way we arrive at a generalized polynomial identity $\sum_i a_{0i}X_1a_{1i}X_2a_{2i}\ldots  a_{n-1i}X_na_{ni}$ (instead of $X_1X_2\ldots X_n$). But using  Theorem \ref{LJerry} this identity can be easily handled.  
Therefore we may assume \eqref{eaf} with both sides nonzero. Repeating the  
 above argument leads us to the situation 
where \cite[Lemmas 6.1.2 and 6.1.4]{BMMb} can be used to prove Martindale's characterization of prime GPI-rings \cite[Theorem 6.1.6]{BMMb}.

It seems that the proof that we outlined is somewhat simpler than the one given in \cite[pp. 216-217]{BMMb}. 
\end{remark}

\section{Prime PI-rings}

We are now in a position to prove the ultimate version of Posner's theorem.

\begin{theorem}\label{main}
If $R$ is a prime PI-ring with extended centroid $C$, then:
\begin{itemize}
\item[(a)] its central closure $R_C$ is a finite-dimensional central simple algebra over $C$;
\item[(b)] every nonzero ideal of $R$  intersects the center $Z$ of $R$ nontrivially;
\item[(c)] $C$ is the field of fractions of $Z$; 
\item[(d)] every element in $R_C$ is of the form $z^{-1} r$ with $0\ne z\in Z$, $r\in R$.
\end{itemize}
\end{theorem}

\begin{proof} (a) Let $U$ be a nonzero ideal of $R_C$. Since $R_C$ is a prime PI-ring (namely, it clearly satisfies the same multilinear identities as $R$), so is $U$.   Let $f = f(X_1,\ldots,X_n)$ be a multilinear polynomial  
identity of $U$ of minimal degree $n$. Write 
$$
f= gX_n + \sum_i g_iX_nh_i
$$
where each $h_i$ is a monomial of degree  $\ge 1$ and with leading coefficient $1$, and $g$ and $g_i$ are  multilinear polynomials. Without loss of generality we may assume that $g\ne 0$. As the degree of $g$ is $n-1$, $g$ is not an identity of $U$.  Pick $u_1,\ldots,u_{n-1}\in U$
so that $u=g(u_1,\ldots,u_{n-1})\ne 0$. We have
$
ux1= ux = \sum v_ixw_i 
$
for some $v_i\in R_C + C$, $w_i\in U$ and all $x\in U$.  Theorem \ref{LJerry} implies that $1$ lies in the $C$-linear span of the $w_i's$. This in particular shows that $1\in R_C$, hence $C\subseteq R_C$, and so $1 \in \sum Cw_i \subseteq U$. Thus $R_C$ is a simple algebra over its center $C$.   Proposition \ref{P1} tells us that it is finite-dimensional.

(b) Let $\varphi$ be a nonzero $C$-linear functional on $R_C$. In view of (iv) in Proposition \ref{P1} there exists 
$T\in M(R_C)$ such that $T(x) = \varphi(x)1$ for all $x\in A$. Let
 $p_i,q_i\in R_C$ be such that $T= \sum_{i=1}^n L_{p_i}R_{q_i}$, $q_1\ne 0$, and $p_1,\ldots,p_n$ are linearly independent (the latter can indeed be required, since otherwise we can pick a maximal linearly independent subset of the $p_i$'s and then rewrite $T$ in an appropriate way).  
Let $J_i$ and $K_i$ be nonzero ideals of $R$ such that $p_iJ_i\subseteq R$ and $K_iq_i\subseteq R$. Now pick any nonzero ideal $I$ of $R$. Then, since $R$ is prime, $I'= (J_1\cap\ldots\cap J_n)I(K_1\cap\ldots \cap K_n)$ is again a nonzero ideal of $R$, and note that $T(I')\subseteq I\cap C$.  Theorem \ref{LJerry} shows that $T(I')\ne 0$, and so $I\cap C\ne 0$.
Since $I\subseteq R$ we actually have $ I\cap C = I\cap Z$.

(c) Let
$\lambda\in C$. Take a nonzero ideal $I$ of $R$ such that $\lambda I\subseteq R$. Picking $0\ne z \in I\cap Z$, we thus have $\lambda z \in R\cap C=Z$. Therefore
$\lambda = z^{-1}z'$ with $z,z'\in Z$.

(d) Use  a common denominator.
\end{proof}

\begin{remark}\label{R3}
A well known result by Rowen \cite{Row} says that (b) holds even for semiprime PI-rings. The standard proof is based on central polynomials. Using this tool - more precisely, we need the existence of a multilinear central polynomial with integer coefficients on a 
 finite-dimensional central simple algebra - one can also derive this more general result by using our approach. Basically one just has to follow Rowen's argument, which, however, can be simplified by omitting  the reduction to the semiprimitive case. Namely,  in view of (a) we may deal with subdirect products of prime rings instead of primitive ones.
 
 Let us point out that we have also used some sort of ``central polynomials" in the proof of (b). However, instead of  usual polynomials we have dealt with a generalized polynomial $\sum_{i=1}^n p_iXq_i$.  As we saw, proving that such a ``polynomial" can have only central values is fairly easy. This cannot be said  for the proof of the existence of the usual central polynomials.
\end{remark}

\begin{remark}\label{R4}
The fact that in our proof we have avoided using the existence of central polynomials on matrix algebras \cite{For, Raz}, makes it possible for us to obtain a new proof of that. Indeed, one just has to use (b) to conclude that the algebra of generic $n\times n$ matrices (which is easily seen to be a prime PI-ring \cite[Corollary 23.52]{Row2}) has a nonzero center; cf. \cite[p.\,324]{Braun}.
\end{remark}

{\bf Acknowledgement.} The author would like to thank Amiram Braun and Louis Rowen for drawing his attention to \cite{Braun}, and to Igor Klep and \v Spela \v Spenko for useful remarks.


\begin{thebibliography}{99}

\bibitem{Bah} Yu. Bahturin, {\em 
Basic structures of modern algebra}, Kluwer, 1993. 

\bibitem{BMMb}
K.\,I. Beidar, W.\,S. Martindale 3rd, A.\,V. Mikhalev, {\em
Rings with generalized identities}, Marcel Dekker, 1996.

\bibitem{Braun}
A. Braun, On Artin's theorem and Azumaya algebras, {\em J. Algebra} {\bf 77} (1982), 223-332.


\bibitem{FIbook} M. Bre\v{s}ar, M.\,A. Chebotar, W.\,S. Martindale 3rd, {\em Functional identities},
Birkh\" auser,  2007.


\bibitem{For}E. Formanek, Central polynomials for matrix rings, {\em J. Algebra} {\bf 23} (1972), 129-132.



\bibitem{Kap}I. Kaplansky, Rings with a polynomial identity, 
{\em Bull. Amer. Math. Soc.} {\bf 54} (1948), 575-580. 


\bibitem{Mart} W.\,S.~Martindale 3rd,
Prime rings satisfying a generalized polynomial identity,
{\em J. Algebra} {\bf 12} (1969), 576-584.

\bibitem{MR} J.\,C.
McConnell, J.\,C. Robson, {\em 
Noncommutative Noetherian rings},
 John Wiley $\&$ Sons, 1987.

\bibitem{Pos}E.\,C. Posner,  Prime rings satisfying a polynomial identity, {\em   Proc. Amer. Math. Soc. } {\bf  11}   (1960), 180-183. 


\bibitem{Raz}Y.\,P. Razmyslov, A certain problem of Kaplansky, {\em Izv. Acad. Nauk SSSR Ser. Mat.} {\bf 37} (1973), 67-74.

\bibitem{Row}L.\,H. Rowen, Some results on the center of a ring with polynomial identity,
{\em Bull. Amer. Math. Soc.} {\bf 79} (1973), 219-223.

\bibitem{Row2}L.\,H. Rowen, {\em  Graduate algebra: noncommutative view},
 American Mathematical Society,  2008.


\end{thebibliography}
\end{document}